\documentclass[12pt]{article}

\usepackage{amsmath,amssymb,amsfonts,theorem,makeidx,latexsym,epsfig}
\usepackage{subfigure}
\usepackage{color}

\newtheorem{defn}{Definition}[section]
\newtheorem{lemma}[defn]{Lemma}
\newtheorem{proposition}[defn]{Proposition}

{\theorembodyfont{\rmfamily}
\newtheorem{ex}[defn]{Example}}
\newtheorem{thm}[defn]{Theorem}
\newtheorem{prop}[defn]{Proposition}
\newtheorem{cor}[defn]{Corollary}

\newcommand{\h}{{\cal H}}
\newcommand{\ltr}{ L^2(\mathbb R) }

\newcommand{\si}{S^{-1}}

\newcommand{\mr}{\mathbb R}

\newcommand{\mz}{\mathbb Z}
\newcommand{\mc}{\mathbb C}

\newcommand{\mts}{ \{E_{mb}T_{na}g \}_{m,n \in \mz}}

\def\bp{{\noindent\bf Proof. \ }}
\def\ep{\hfill$\square$\par\bigskip}

\def\bqs{\begin{equation}}
\def\eqs{\tag*{$\square$}\end{equation}\par\bigskip}

\def\la{\langle}
\def\ra{\rangle}

\def\nl{\left|\left|}
\def\nr{\right|\right|}
\def\span{\overline{\text{span}}}

\def\bee{\begin{eqnarray}}
\def\ene{\end{eqnarray}}

\def\bes{\begin{eqnarray*}}
\def\ens{\end{eqnarray*}}

\def\bei{\begin{itemize}}
\def\eni{\end{itemize}}

\def\bt{\begin{thm}}
\def\et{\end{thm}}

\def\bc{\begin{cor}}
\def\ec{\end{cor}}

\def\bpr{\begin{prop}}
\def\epr{\end{prop}}

\def\bl{\begin{lemma}}
\def\el{\end{lemma}}

\def\bd{\begin{defn}}
\def\ed{\end{defn}}

\def\bex{\begin{ex}}
\def\enx{\end{ex}}

\def\bfi{\begin{fig}}
\def\efi{\end{fig}}

\newcommand{\nft}{ || f||^2}

\def\eti{\{e_i\}_{i\in I}}

\def\hti{\{h_i\}_{i\in I}}

\def\fti{\{f_i\}_{i\in I}}

\def\cti{\{c_i\}_{i\in I}}

\def\otj{\{\omega_j\}_{j\in I}}

\def\sui{\sum_{i\in I}}

\newcommand{\seqgri[1]}{\{#1_i\}_{i\in I}}

\newcommand{\seqgrj[1]}{\{#1_j\}_{j\in I}}

\newcommand{\sumgr[1]}{\sum_{#1=1}^\infty}
\def\newin {\,\kern-0.4em\in\kern-0.15em}
\def\newsubset {\kern-0.2em\subset\kern-0.2em}

\def\sumgri{\sum_{i\in I}}
\def\sumgrj{\sum_{j\in I}}

\def\<{\langle}
\def\>{\rangle}

{\nopagebreak\par\noindent\hbox to \hsize{\hrulefill}\par\noindent}

\title{On various R-duals and the duality principle}

\date{\today}

\author{Diana Stoeva, Ole Christensen}

\begin{document}

\maketitle

\begin{abstract} The duality principle states that a Gabor system is a frame if and only if the corresponding adjoint Gabor system is a Riesz sequence. In general Hilbert spaces and without the assumption of any particular structure, Casazza, Kutyniok and Lammers have introduced the so-called R-duals that also lead to a characterization of frames in terms of associated Riesz sequences; however, it is still an open question whether this abstract theory is a generalization of the duality principle.
In this paper we prove that a modified version of the R-duals leads to a generalization of the duality principle that keeps all the attractive properties of the R-duals. In order to provide extra insight into the relations between a given sequence and its R-duals, we characterize
all the types of R-duals that are available in the literature for the special case where the
underlying sequence is a Riesz basis.
\end{abstract}


\vspace{.05in}

\section{Introduction}
A countable collection of vectors $\fti$ in a separable Hilbert space $\h$ is a {\it frame for $\h$ with (frame) bounds $A,B$} if 
$A$ and $B$ are strictly positive constants and the inequalities
\bee \label{50329a} A\, \nft \le \sum_{i\in I} | \la f, f_i\ra|^2 \le B\, \nft\ene
hold for all $f\in \h.$ 
 Frames play an increasing role in analysis and applications, mainly due to the fact that frames yield expansions of the elements in the Hilbert space of a
similar type as the one that is known for orthonormal bases. In fact, if $\fti$ is a frame for $\h,$ the frame operator $S: \h \to \h, \, Sf:=  \sum_{i\in I} \la f, f_i\ra f_i$, is known to be invertible, and
\begin{equation} \label{dual}
 f= \sum_{i\in I} \la f, \si f_i \ra f_i,
 \end{equation} for all $f\in\h$.
It is clear that it might be a nontrivial matter to verify the two inequalities in \eqref{50329a}.  For so-called Gabor systems in $\ltr$
(see the description below),
the
{\it duality principle}  \cite{DLL,Jan5,RoSh1} states that
the frame condition is equivalent to a Riesz basis condition on an associated sequence
(the adjoint Gabor system), see Theorem \ref{1505d1}; this leads to a method to check the frame condition for Gabor systems
in an (at least conceptually) easier way.  In an attempt to extend this to general
sequences in arbitrary Hilbert spaces, Casazza, Kutyniok and Lammers introduced the
{\it R-duals} in the paper \cite{CKL}. The R-duals also yield a method for checking the
frame condition for a sequence of vectors by checking the Riesz basis condition for a related sequence.

At present it is not known whether the theory for R-duals yields a generalization of the duality principle.
In \cite{sc} the authors introduced certain variations of the R-duals  (see Definition \ref{defrduals}) and showed that
{\it R-duals of type II} cover the duality principle for integer-oversampled Gabor systems leaving open the general case,
while {\it R-duals of type III} generalize the duality principle  and keep some of the attractive properties of the R-duals, but not all.
In the current paper we show that R-duals of type II in fact do not generalize the
duality principle for arbitrary Gabor frames. This brings the attention to the R-duals of type III and
we determine a sub-class of the R-duals of type III, which possesses the missing properties.
We also provide further insight into the various R-duals by providing characterizations in the special case
where the given frame is a basis.

In the rest of this introduction we state the key definitions and results from the
literature concerning the R-duals. In Section \ref{50329b} we introduce the modified
R-duals of type III and prove that they generalize the duality principle and keep the
main properties known from the Gabor case.  The special case of Riesz bases is analysed
in Section \ref{rbcase}.
We refer to the monographs
\cite{Cbook2,G2,Heil} for detailed treatments of frames and further references.

 For a sequence
$\fti$ satisfying at least the upper frame condition,
the {\it analysis operator} is defined by
\bes U: \h \to \ell^2(I), \, Uf:= \{\la f, f_i\ra\}_{i\in I}\ens
and the {\it synthesis operator is}
\bes T: \ell^2(I) \to \h, \, T \cti = U^*\cti = \sum_{i\in I} c_i f_i.\ens
A sequence $\seqgri[f]$ with elements from $\h$ is called a {\it frame sequence in $\h$} if the inequalities in \eqref{50329a}
hold for all $f\in \span\fti$.
Given a frame sequence $\seqgri[f]$ in $\h$, its {\it frame operator} is $S: \span \fti \to \span \fti, \,
Sf:=\sum_{i\in I}\<f,f_i\>f_i.$
 In this case $S$ is a bijection on $\span \fti$
and
the sequence $\{\si f_i\}_{i\in I}$ is  a frame for $\span\fti$ satisfying the representation formula (\ref{dual}) for all $f\in\span \fti$;
 $\{\si f_i\}_{i\in I}$ is called
the {\it  canonical dual frame of $\seqgri[f]$}
and will be denoted by $\seqgri[\widetilde{f}]$.
A sequence $\fti$ with elements from $\h$ is called a {\it Riesz sequence in $\h$ with bounds $A,B$} 
if $0<A\leq B<\infty$ and 
\bes A \, \sum_i |c_i|^2 \le \nl \sum_i c_i f_i \nr^2 \le B\, \sum_i |c_i|^2\ens for all finite sequences $\{c_i\}_{i\in I}.$ A Riesz sequence is called a {\it Riesz basis for $\h$} if
$\h= \span \fti.$
Recall that if $\fti$ is a frame for $\h$, then $\|S\|$ is the optimal upper frame bound and $\|S^{-1}\|^{-1}$ is the optimal lower frame bound.

Let us now collect the various definitions of R-duals that are available in the literature.
The definition by Casazza \& al. corresponds to what we call R-duals of type I.

\begin{defn}\label{defrduals}
Let $\eti$ and $\hti$ be sequences with elements in $\h$ and let $\fti$ be a sequence in $\h$ for which $\label{14a} \sui
|\la f_i,e_j\ra|^2 < \infty, \  \forall j\in I.$
\bei
 \item[{\rm (i)}] {\rm \cite{CKL}}
When $\eti$ and $\hti$ are orthonormal bases
for $\h,$ the {\bf R-dual of type I} of $\seqgri[f]$  with respect to {\rm ($\seqgri[e]$, $\seqgri[h]$)} is the sequence  $\seqgrj[\omega]$ given by
\begin{equation}\label{grdual1}
 \omega_j=\sumgri \<  f_i,  e_j\>  h_i, \ j\in I.
 \end{equation}
\item[{\rm (ii)}] {\rm \cite{sc}} Let $\seqgri[e]$ and $\seqgri[h]$ be orthonormal bases for $\h.$ If $\seqgri[f]$ is a frame for $\h$ with frame operator $S$, the {\bf R-dual of type II} of $\seqgri[f]$ with respect to {\rm ($\seqgri[e]$, $\seqgri[h]$)} is the sequence  $\seqgrj[\omega]$ given by
\begin{equation}\label{grdual}
 \omega_j=\sumgri \<  f_i, S^{-1/2} e_j\> S^{1/2} h_i, \ j\in I.
 \end{equation}
\item[{\rm (iii)}] {\rm \cite{sc}} Let $\seqgri[e]$ and $\seqgri[h]$ be orthonormal bases for $\h$.
If $\seqgri[f]$ is a frame sequence in $\h$ with frame operator $S$ and $Q:\h\to\h$ is a bounded bijective
operator with  $\|Q\|\leq \sqrt{||S||}$ and  $\|Q^{-1}\|\leq \sqrt{|| \si ||}$, the {\bf R-dual of type III} of $\seqgri[f]$ with respect to the triplet {\rm($\seqgri[e]$, $\seqgri[h]$, $Q$)} is the sequence  $\seqgrj[\omega]$ defined by
\begin{equation}\label{rfdualGab}
\omega_j:=\sumgr[i] \< S^{-1/2}f_i, e_j\> Qh_i \ j\in I.
\end{equation}
\item[{\rm (iv)}]
  {\rm \cite{XZ}} When $\seqgri[e]$ and $\seqgri[h]$ are Riesz bases for $\h$, the {\bf R-dual of type IV} of $\seqgri[f]$ with respect to {\rm($\seqgri[e]$, $\seqgri[h]$)} is the sequence  $\seqgrj[\omega]$ given by
\begin{equation}\label{grdualp}
 \omega_j=\sumgri \<  f_i, e_j\> h_i, \ j\in I.
 \end{equation}
\eni
\end{defn}

R-duals of type I are interesting because they form Riesz sequences if and only
if the given sequence $\fti$ is a  frame for $\h:$

\bt \label{207b} {\rm \cite{CKL}} Let $\fti$ be a sequence in $\h$ and let
$\otj$ be an R-dual of type I of $\fti$. Then the following hold:
\begin{itemize}
\item[{\rm (i)}]  $\fti$ is a frame for $\h$ with bounds $A,B$ if and
only if $\otj$ is a Riesz sequence in $\h$  with bounds $A,B.$
\item[{\rm (ii)}]  $\seqgrj[\omega]$  is a Riesz basis for $\h$ if and only if $\seqgri[f]$ is a Riesz basis for $\h$.
\end{itemize}
 \et

In the literature, several characterizations of the various types of R-duals are
formulated in terms of the condition
\bee \label{223f} {\rm dim (ker} \, T)= {\rm dim ((span}\otj)^\perp),\ene relating the
synthesis operator $T$ for the sequence $\fti$ to the sequence $\otj.$ We
collect these results here:

\bl\label{impcond}
Let $\seqgri[f]$ be a  frame for $\h$ and let $\otj$ be a Riesz sequence in $\h$.
Denote the synthesis operator for $\fti$ by  $T$, the frame operator of  $\fti$ by  $S$, and the frame operator of $\otj$ by $S_\Omega$. The following statements hold.
\begin{itemize}
\item[{\rm (i)}] {\rm \cite{CKL}} If  $\otj$ is an R-dual of type I of $\fti$ with respect to {\rm ($\seqgri[e]$, $\seqgri[h]$)}, then
$g\in ({\rm span} \otj)^\perp$ if and only if $\{\<h_i,g\>\}_{i\in I}\in {\rm ker} \, T;$
in the affirmative case  (\ref{223f}) holds.
\item[{\rm (ii)}] {\rm \cite{sc}} If $\seqgri[f]$ is tight and  $\otj$ is tight with the same 
bound, then $\otj$ is an R-dual of type I of $\fti$ if and only if (\ref{223f}) holds.
\item[{\rm (iii)}] {\rm \cite{cz}}  $\otj$ is an R-dual of type I of $\fti$ if and only if (\ref{223f}) holds and there exists an antiunitary transformation $\Lambda: \h \to \overline{{\rm span}}\otj$ so that $S_\Omega=\Lambda S \Lambda^{-1}$.
\item[{\rm (iv)}] {\rm \cite{sc}}
 $\otj$ is an R-dual of type II of $\fti$ if and only if (\ref{223f}) holds, $\{S^{-1/2}\omega_j\}_{j\in I}$ is orthonormal,  and the frame bounds of $\seqgri[f]$ are also 
 bounds of $\otj$.
\item[{\rm (v)}] {\rm \cite{sc}}
 $\otj$ is an R-dual of type III of $\fti$ if and only if (\ref{223f}) holds
 and the frame bounds of $\seqgri[f]$ are also 
 bounds of $\otj$.

\end{itemize}

\el

In the discussion of the duality principle we need the definition of a Gabor system.
Consider the Hilbert space $L^2(\mr).$ For $p,q\in\mr$, let
$T_p: L^2(\mr)\to L^2(\mr)$ denote the translation operator  $(T_p g) (x) = g(x-p)$, and  $E_q: L^2(\mr)\to L^2(\mr)$ the modulation operator $(E_q g)(x) = e^{2\pi  {\rm i} q x} g(x).$
Given parameters $a>0, b>0$ and $g\in \ltr,$ the associated
{\it Gabor system}  is the sequence  $\{E_{mb} T_{na} g\}_{m,n\in\mz};$ the
{\it adjoint Gabor system}  is the sequence $\{E_{m/a} T_{n/b} g\}_{m,n\in\mz}.$

The duality principle, due to Janssen \cite{Jan5}, Daubechies, Landau,
and Landau \cite{DLL}, and Ron and Shen \cite{RoSh1}, states the following:

\begin{thm} \label{1505d1} {\rm \cite{DLL, Jan5, RoSh1}} Let $g\in \ltr$ and $a,b>0$ be given.
Then the Gabor system $\mts$ is a frame for $L^2(\mr)$ with bounds $A,B$ if and only if $\{\frac1{\sqrt{ab}} \,
E_{m/a}T_{n/b}g\}_{m,n\in \mz}$ is a Riesz sequence with bounds $A,B.$
\end{thm}
It is well known that if a Gabor system $\mts$  is a Riesz basis, then $ab=1;$ via the
duality principle this implies that the Gabor system $\mts$ is a Riesz basis for $L^2(\mr)$ if and only if $\{\frac1{\sqrt{ab}} \,
E_{m/a}T_{n/b}g\}_{m,n\in \mz}$ is a Riesz basis for $L^2(\mr)$. Thus, the relation between the
Gabor system $\mts$ and the system $\{\frac1{\sqrt{ab}} \,
E_{m/a}T_{n/b}g\}_{m,n\in \mz}$ corresponds exactly to the relation between the sequence $\fti$
and its R-duals, see Theorem \ref{207b}. However, it is still not known whether this is a coincidence, or the R-duals of type I actually generalize the duality principle. That is, given a Gabor frame $\mts$ for $L^2(\mr)$,  we do not know whether the Gabor system
$\{ \frac1{\sqrt{ab}} E_{m/a}T_{n/b}g\}_{m,n\in \mz}$ always can be realized as an R-dual of type I of $\mts$ (by \cite{CKL}, the answer is affirmative for  tight Gabor frames and Gabor Riesz bases).
This is the motivation for the introduction of the other types of R-duals in \cite{sc}. In
fact,
in \cite{sc} it was shown that the R-duals of type III generalize the duality principle for
all Gabor systems.
However, R-duals of type III do not enjoy all of the attractive properties of the duality
principle, so it is natural and necessary to search for subclasses which are in closer correspondence with
the duality principle. In the present paper we determine a relevant subclass of the R-duals of type III which both extends the duality principle and has the desired properties as in Theorem \ref{207b}.

 Note that the duality principle and the R-duals by Casazza et al.
 have trigged a lot of research activity; we refer to
 the papers \cite{CXZ,cz, H, FS, fhs, XZ}.

\section{R-duals of type  II} \label{gabsection}

In this section we solve one of the remaining problems
in \cite{sc}, by showing that the  R-duals of type II
do not generalize the duality principle. This will motivate the analysis in the rest of the paper, where
we focus on a subclass of the R-duals of type III having particular properties.

\begin{ex} \label{rd2counterex}
Consider the B-spline $B_2$ determined by

   $$B_2(x)=\left\{
 \begin{array}{ll}
 1+x,&  x\in[-1, 0];\\
 1-x, & x\in[0, 1];\\
 0, & x\notin [-1,1].
   \end{array}
 \right.
 $$
   Let $a=1$ and $b=\frac{2}{5}$. 
 By \cite[Corollary 11.7.1]{Cbook2}, $\{E_{mb}T_{na} B_2\}_{m,n\in \mz}$ is a Gabor frame for $L^2(\mr)$.
We will now prove that the Riesz sequence  $\{\frac{1}{\sqrt{ab}} E_{m/a}T_{n/b} B_2\}_{m,n\in \mz}$
 is not an R-dual of type II of $\mts$.

Denote the frame operator of $\{E_{mb}T_{na} B_2\}_{m,n\in \mz}$  by $S$.
By \cite[Corollary 11.4.5]{Cbook2}, 
 $Sf=\frac{G}{b}f \ \mbox{ and } \ S^{-1}f=\frac{b}{G}f, \ f\in L^2(\mr),$ where
$G(x)=\sum_{k\in\mz} |B_2(x-ka)|^2, \ x\in\mr.$
Denote $\omega_{m,n}:=\frac{1}{\sqrt{ab}} E_{m/a}T_{n/b} B_2$, $m,n\in \mz$. 
We will show that $\{ S^{-1/2} \omega_{m,n}\}_{m,n\in \mz}$ is not an orthonormal sequence, which by Lemma \ref{impcond}(iv) will imply that $\{\omega_{m,n}\}_{m,n\in \mz}$ is not an R-dual of type II of $\mts$. For $m=0$ and $n=1$, we have
 \begin{eqnarray*}
\< S^{-1/2} \omega_{0,1},  S^{-1/2} \omega_{0,1}\> & =& \< S^{-1} \omega_{0,1},  \omega_{0,1}\> 
= \int_{\mr} \frac{b}{G(x)} \omega_{0,1}(x)  \overline{\omega_{0,1}(x)}\, {\rm d x}  \\
 & =&  \int_{\mr} \frac{b}{G(x)} |\omega_{0,1}(x)|^2 \, {\rm d x}.
\end{eqnarray*} Since
  $$
  \omega_{0,1}(x) = \frac{1}{\sqrt{ab}} E_0 T_{1/b} B_2 (x) 
= \frac{1}{\sqrt{b}} B_2 (x-\frac{5}{2}) =
\frac{1}{\sqrt{b}}\left\{
 \begin{array}{ll}
 x-\frac{3}{2},&  x\in[\frac{3}{2}, \frac{5}{2}],\\
 \frac{7}{2}-x, & x\in[\frac{5}{2}, \frac{7}{2}],\\
 0, & x\notin [\frac32, \frac72],
   \end{array}
 \right.
 $$
and for $x\in[\frac{3}{2},\frac{7}{2}]$,
 $$G (x)=\sum_{k\in\mz} |B_2(x-k)|^2 =
 \left\{
 \begin{array}{ll}
 (x-1)^2 + (2-x)^2,&  x\in[\frac{3}{2}, 2];\\
 (x-2)^2 + (3-x)^2,&  x\in[2, 3]; \\
 (x-3)^2 + (4-x)^2,&  x\in[3, \frac{7}{2}], \\
   \end{array}
 \right.
 $$
it follows that

 \begin{eqnarray*}
 \< S^{-1/2} \omega_{0,1},  S^{-1/2} \omega_{0,1}\>   =
\int_{\mr} \frac{1}{G(x)} |B_2 (x-5/2)|^2 \, {\rm d x} = 1+\frac{\pi}{4} - \ln 2 \neq 1. 
  \end{eqnarray*}
Therefore, $\{\omega_{m,n}\}_{m,n\in\mz}$ is not an R-dual of type II of $\mts$. \ep \enx

\section{R-duals of type III} \label{50329b}

The key motivation behind the definition of R-duals of type III is that they generalize the duality principle \cite{sc}, in the sense that whenever $\mts$ is a frame for $L^2(\mr)$, the system $\{\frac1{\sqrt{ab}} \,
E_{m/a}T_{n/b}g\}_{m,n\in \mz}$ can be realized as an R-dual of type III of $\mts$. However, not all
R-duals of type III have exactly the same properties as encountered in the duality principle.
For example, for a frame $\fti$ with frame operator $S,$ the optimal frame bounds are
$\frac{1}{\|S^{-1}\|}$, $ \|S\|;$ these numbers are also bounds for the R-duals of type III, but
not necessarily the optimal bounds. This calls for an identification of a subclass of the
R-duals of type III with properties that better match what we know from the duality principle.

As starting point
we will now determine  conditions on the operator $Q$ in \eqref{rfdualGab} which are necessary and sufficient for an R-dual of type III of $\fti$ to keep the optimal bounds of $\fti$.

\begin{prop}  \label{rdualopt}
Let $\fti$ be a frame for $\h$ (resp. Riesz sequence in $\h$) with frame operator $S$ and analysis operator $U,$ both considered as operators on $\span \fti$, and let $\otj$ be an R-dual of type III of $\fti$
with respect to the triplet {\rm($\seqgri[e]$, $\seqgri[h]$, $Q$)}.
Then the following are equivalent:
\bei \item[(i)] The Riesz sequence (resp. the frame) $\otj$
has the same optimal bounds as $\fti.$
\item[(ii)] The operator $Q$ has the property 
\begin{equation}  \label{eqstar}
\left\{ \begin{array}{l}
\frac{1}{\sqrt{\|S^{-1}\|}} \|\seqgri[d]\|_{\ell^2} \leq \|Q(\sum_{i\in I} d_i h_i) \| \leq \sqrt{\|S\|} \|\seqgri[d]\|_{\ell^2}
\\
\mbox{ whenever } \{\overline{d_i}\}_{i\in I}\in {\mathcal R}(U), \mbox{with optimality of the bounds} \frac{1}{\sqrt{\|S^{-1}\|}}, \sqrt{\|S\|}.
\end{array}
\right.
\end{equation}
\eni
\end{prop}

\bp Notice that when $\fti$ is a frame for $\h$ (resp. Riesz sequence),
\cite[Prop. 4.3]{sc} shows that $\otj$ is a Riesz sequence (resp. frame for $\h$), with bounds $1/\|S^{-1}\|$, 
$\|S\|$. We first consider the case where $\fti$ is assumed to be a frame for $\h$.

\noindent (i) $ \Rightarrow$ (ii) Assume that $\frac{1}{\|S^{-1}\|}$ and $\|S\|$ are the optimal bounds of $\otj$. We will prove that (\ref{eqstar}) holds.
 Since
$\|Q\|\leq \sqrt{||S||}$ and  $\|Q^{-1}\|\leq \sqrt{|| \si ||}$, it follows that the inequalities
of (\ref{eqstar})
 hold for all $\{d_i\}_{i\in I}\in\ell^2$ and in particular whenever $\{\overline{d_i}\}_{i\in I}\in {\mathcal R}(U)$; it remains to prove the optimality of the bounds of (\ref{eqstar}).
Assume that there exists $B_1< \sqrt{\|S\|}$ so that
$$ \|Q(\sum_{i\in I} d_i h_i) \| \leq B_1 \|\seqgri[d]\|_{\ell^2}
\mbox{ whenever } \{\overline{d_i}\}_{i\in I}\in {\mathcal R}(U).$$
Then for every finite scalar sequence $\{c_j\}$, taking
\bes \{\overline{d_i}\}_{i\in I} =
 \{\<\sum_j \overline{c_j} e_j,  S^{-1/2}f_i\>\}_{i\in I}
 =\{\<S^{-1/2}(\sum_j \overline{c_j} e_j),  f_i\>\}_{i\in I}
 \in R(U),\ens
we obtain
 \begin{eqnarray*}
 \|\sum_j c_j \omega_j \| &= &
  \|Q(\sum_{i\in I} \< S^{-1/2}f_i, \sum_j \overline{c_j} e_j\> h_i) \| = \|Q(\sum_{i\in I} d_i h_i) \|
  \leq
  B_1 \|\seqgri[d]\|_{\ell^2}
 \\
 & =&
  B_1 \|\{\<\sum_j \overline{c_j} e_j, S^{-1/2} f_i\>\}_{i\in I}\|_{\ell^2}
= B_1 \|\{c_j\}\|_{\ell^2}. \end{eqnarray*}
This implies that $B_1^2$ is an upper bound of the Riesz sequence $\otj$, which contradicts  the assumptions. Therefore, $\sqrt{\|S\|}$ is the optimal upper bound  in (\ref{eqstar}). In a similar way, it follows that $1/\sqrt{\|S^{-1}\|}$ is the optimal lower bound in (\ref{eqstar}).

\noindent (ii) $ \Rightarrow$ (i) Now assume that (\ref{eqstar}) holds; we will prove that $1/\|S^{-1}\|$ and $\|S\|$ are the optimal bounds of the Riesz sequence $\otj$.
As already mentioned, $1/\|S^{-1}\|$ and $\|S\|$ are bounds of the Riesz sequence $\otj$, so it remains to prove their optimality.
 Assume that the optimal upper bound of  $\otj$ is $B_2$ with  $B_2< \|S\|$.
 Then for every finite scalar sequence $\{c_j\}$, we have
 \begin{eqnarray*}
   \|Q(\sum_{i\in I} \< S^{-1/2}f_i, \sum_j c_j e_j\> h_i) \|
   &= &
 \|\sum_j \overline{c_j} \omega_j \|
  \leq
  \sqrt{B_2}\, \|\{\overline{c_j}\}\|_{\ell^2}
   = \sqrt{B_2}\, \|\{c_j\}\|_{\ell^2}
 \\
 & =&
  \sqrt{B_2}\,\|\{\<S^{-1/2}f_i, \sum_j c_j e_j\>\}_{i\in I}\|_{\ell^2}. \end{eqnarray*}
Since the set of (finite) linear combinations $\sum_j c_j e_j$ is dense in $\h$, it follows that
  \begin{eqnarray*}
   \|Q(\sum_{i\in I} \< S^{-1/2}f_i, y\> h_i) \|
  \leq
   \sqrt{B_2} \,  \|\{\<S^{-1/2}f_i, y\>\}_{i\in I}\|_{\ell^2}   \end{eqnarray*}
for every $y\in\h$.
Since $S^{-1/2}$ is bijective and self-adjoint, it follows that   $
   \|Q(\sum_{i\in I} \< f_i, u\> h_i) \|
  \leq
   \sqrt{B_2} \, \|\{\< f_i, u\>\}_{i\in I}\|  $  for all $u\!\in\!\h$, which contradicts (\ref{eqstar}).

Now let $\fti$ be a Riesz sequence in $\h$.
In this case  $R(U)=\ell^2$; since the optimal bounds in the inequalities
$C ||x|| \le ||Qx|| \le D ||x||, x\in \h,$ are $C= ||Q^{-1}||^{-1}, D=||Q||,$ the condition (\ref{eqstar}) means precisely that $\|Q\|= \sqrt{||S||}$ and  $\|Q^{-1}\|= \sqrt{|| \si ||}$.
An argument as in the proof of \cite[Prop. 4.3(ii)]{sc} shows that $\otj$ is a frame with optimal bounds $\frac{1}{\|Q^{-1}\|^2}$, $\|Q\|^2$.
Therefore, $\otj$ has optimal bounds $\frac{1}{\|S^{-1}\|}$ and $\|S\|$ if and only if  $\|Q^{-1}\|= \sqrt{|| \si ||}$ and $\|Q\|= \sqrt{||S||},$ i.e., if and only if
(\ref{eqstar}) holds.
\ep

Proposition \ref{rdualopt} identifies the relevant subclass of the R-duals of type III.  Fortunately
we can now show that this class keeps the essential property of generalizing the duality principle:

\begin{thm} \label{GaborIII}
Let $\mts$ be a Gabor frame for  $L^2(\mr)$.
Then \\ $\{\frac{1}{\sqrt{ab}} E_{m/a}T_{n/b} g\}_{m,n\in \mz}$ can be realized as an
R-dual of type III of $\mts$ with an operator $Q$ having the property (\ref{eqstar}).
 \end{thm}
\bp By the duality principle (Theorem \ref{1505d1}), $\{\frac{1}{\sqrt{ab}} E_{m/a}T_{n/b} g\}_{m,n\in \mz}$ is a Riesz sequence which has the same optimal bounds as  $\mts$.
By \cite[Corollary 4.5]{sc},  $\{\frac{1}{\sqrt{ab}} E_{m/a}T_{n/b} g\}_{m,n\in \mz}$ can be realized as an
R-dual of type III of $\mts$ with respect to  some orthonormal bases $\seqgri[e]$, $\seqgri[h]$, and
an appropriate operator $Q$.
Now by Proposition \ref{rdualopt},
this operator $Q$
 must satisfy the property  (\ref{eqstar}).
\ep

Furthermore, the following result (which is an immediate consequence of Proposition \ref{rdualopt} and \cite[Prop. 4.3]{sc})
shows that the class of
R-duals of type III having the property (\ref{eqstar}) provides us with exactly the same frame bounds as the given frame:

\begin{thm} \label{cordualopt} Let $\seqgri[f]$ be a frame  sequence in $\h$ and let $\seqgri[\omega]$ be an R-dual of $\fti$ of type III  with the property (\ref{eqstar}). Then  the following holds.
\bei
\item[(i)] $\seqgri[f]$ is a frame for $\h$ with bounds $A, B$ if and only if $\otj$ is a Riesz sequence with bounds $A, B$;
\item[(ii)] $\seqgri[f]$ is a Riesz sequence with bounds $A, B$ if and only if $\otj$ is a frame for $\h$ with bounds $A, B$;
\item[(iii)] $\seqgri[f]$ is a Riesz basis for $\h$
with bounds $A, B$
if and only if $\otj$ is a Riesz basis for $\h$
 with bounds $A, B$.
 \eni
\end{thm}

We have now identified the correct subclass of the R-duals of type III. It has a compact
characterization in terms of the condition \eqref{223f}:

\begin{thm} \label{symm2}
Let $\seqgri[f]$ be a frame for $\h$ and let $\seqgrj[\omega]$ be a Riesz sequence in $\h$ with the same optimal bounds.
Then (\ref{223f}) holds if and only if
  $\otj$ is an R-dual  of type III of $\fti$  having the property (\ref{eqstar}).
\end{thm}

\bp First assume that $\otj$ is an R-dual of type III of $\fti$ having the property (\ref{eqstar}).
Then  \cite[Theorem 4.4]{sc} implies that (\ref{223f})  holds.

Conversely, assume that (\ref{223f})  holds.
By \cite[Theorem 4.4(ii)]{sc}, $\otj$ is an R-dual of type III of $\fti$
with respect to 
an appropriate triplet  {\rm($\seqgri[e]$, $\seqgri[h]$, $Q$)}.
By  Proposition \ref{rdualopt}, 
the operator $Q$ must satisfy (\ref{eqstar}).
\ep

The next result relates the class of R-duals of type III having the property (\ref{eqstar})
with the R-duals of type I and III, respectively.

\begin{proposition} \label{213a2iii}
Let $\seqgri[f]$ be a  frame for $\h$. Then the following holds.
\begin{itemize}
\item[{\rm (i)}] The class of type I duals of $\fti$ is contained in the class of type III duals having the property (\ref{eqstar}).
\item[{\rm (ii)}] When $\fti$ is tight, the classes mentioned in (i) coincide.
\item[{\rm (iii)}] When $\fti$ is not tight, the class of R-duals of type III having the property (\ref{eqstar}) is a strict subset of the class of R-duals of type III.
\end{itemize}
\end{proposition}

\bp (i) Let $\otj$ be an R-dual of type I of $\fti$. By Theorem \ref{207b}, $\otj$ is a Riesz sequence in $\h$ and
the optimal bounds of $\otj$ are the same as the optimal ones of $\fti$. By \cite[Theorem 4.4(iii)]{sc},
$\otj$ can be written as an R-dual of type III of $\fti$ with respect to some triplet ($\seqgri[e], \seqgri[h], Q$) and by Proposition \ref{rdualopt}, the property (\ref{eqstar}) must hold.

(ii) Assume that $\fti$ is tight.
Then the classes of R-duals of type I and type III of  $\fti$ coincide \cite{sc}.
Now the statement follows from (i).

(iii)
Assume that $\fti$ is not tight and let $A$ and $B$  denote the optimal bounds of $\fti$, $A<B$. Take any constant $C\in (\sqrt{A}, \sqrt{B})$ and let $Q:= C \, {\rm Id}_{\h}$. Let $\otj$ be an R-dual of type III with respect to some orthonormal bases $\eti$, $\seqgri[h]$ and the operator $Q$. Then $\{C^{-1}\omega_j\}_{j\in I}$ is an R-dual of type I of $\{S^{-1/2}f_i\}_{i\in I}$, which by
Theorem \ref{207b}
implies that $\{C^{-1}\omega_j\}_{j\in I}$ is an orthonormal sequence. Therefore, $\otj$ is a tight Riesz sequence with bound $C^2\in (A,B)$, which by Theorem \ref{cordualopt}(i) implies that $\otj$ can not be written as an R-dual of type III with property (\ref{eqstar}). 
\ep

In  \cite{sc} we have proved that canonical dual frames lead to biorthogonality of appropriately determined R-duals of type III. 
 Here we provide further insight in the relations considering the converse situation, namely,
 biortogonality of appropriate R-duals of type III leading to canonical dual frames.

 \begin{prop} \label{biorth}
 Let $\seqgri[f]$ be a frame for $\h$ with frame operator $S$ and analysis operator $U$. 
The following holds.
\begin{itemize}
\item[{\rm (i)}] If $\seqgrj[\omega]$ is an R-dual of type III  of $\seqgri[f]$, then the  biorthogonal sequence of $\seqgrj[\omega]$ in
$\span {\seqgrj[\omega]}$
 is an R-dual of type III  of $\{\widetilde{f_i}\}_{i\in I}$.
 \item[{\rm (ii)}] If $\seqgrj[\omega]$ is an R-dual of type III  of $\seqgri[f]$ with respect to  {\rm($\seqgri[e]$, $\seqgri[h]$, $Q$)} having the property (\ref{eqstar}), then
 the  biorthogonal sequence of $\seqgrj[\omega]$ in
$\span {\seqgrj[\omega]}$
 is an R-dual of type III  of $\{\widetilde{f_i}\}_{i\in I}$ with respect to a triplet {\rm($\seqgri[e]$, $\seqgri[z]$, $V$)}such that
 \begin{equation*}
\left\{ \begin{array}{l}
\frac{1}{\sqrt{\|S\|}} \|\seqgri[d]\|_{\ell^2} \leq \|V(\sum_{i\in I} d_i z_i) \| \leq \sqrt{\|S^{-1}\|} \|\seqgri[d]\|_{\ell^2}
\\
\mbox{ whenever } \{\overline{d_i}\}_{i\in I}\in {\mathcal R}(U), \mbox{with optimality of the bounds}.
\end{array}
\right.
\end{equation*}
 \end{itemize}
 \end{prop}

\bp
Let $\seqgrj[\omega]$ be an R-dual of type III of $\seqgri[f]$
with respect to the triplet {\rm($\seqgri[e]$, $\seqgri[h]$, $Q$)}.
The only biorthogonal sequence of $\seqgrj[\omega]$ in
$\span {\seqgrj[\omega]}$ is the canonical dual frame $\{\widetilde{\omega_j}\}_{j\in I}$ of $\otj$.
Let 
$S_{\widetilde{f}}$, $S_\omega$, and $S_{\widetilde{\omega}}$ denote the frame operators of 
$\{\widetilde{f_i}\}_{i\in I}$, $\otj,$ and $ \{\widetilde{\omega_j}\}_{j\in I}$ respectively.

By \cite[Prop. 4.3(i)]{sc} and some arguments from the proof of \cite[Theor. 4.4(ii)]{sc}, it follows that the 1-tight Riesz sequence $\{S_\omega^{-1/2} \omega_j\}_{j\in I}$ is an R-dual of type I of the 1-tight frame $\{S^{-1/2} f_i\}_{i\in I}$.
Now the proof of \cite[Prop. 1.6]{sc} implies that for the fixed
orthonormal basis $\seqgri[e]$, there exists an orthonormal basis $\seqgri[z]$ of $\h$ so that
 $\{S_\omega^{-1/2} \omega_j\}_{j\in I}$ is an R-dual of type I of  $\{S^{-1/2} f_i\}_{i\in I}$ with respect to $(\seqgri[e], \seqgri[z])$.
Then for every $i\in I$, one can write
\begin{eqnarray*}
{\widetilde{f}}_i &=& 
\sum_{j\in I}  \<S^{-1/2}f_i, e_j\> S^{-1/2} e_j
=\sum_{j\in I} \left\< \sum_{k\in I} \<S^{-1/2}f_k, e_j\>z_k, z_i\right\> S^{-1/2}e_j\\
&=&\sum_{j\in I} \< S^{-1/2}_{\omega} {\omega_j}, z_i\> S^{-1/2}e_j
= \sum_{j\in I} \< S^{-1/2}_{\widetilde{\omega}} \widetilde{\omega_j}, z_i\> S_{\widetilde{f}}^{1/2}e_j,
\end{eqnarray*}
which implies that
\begin{equation*}
S^{-1/2}_{\widetilde{\omega}} \widetilde{\omega_j} =
\sum_{i\in I} \< S_{\widetilde{f}}^{-1/2} \widetilde{f}_i, e_j\>  z_i, \ j\in I.
\end{equation*}
By \cite[Lemma 1.3]{sc}, there exists a bounded bijective
extension $V$ 
of $S^{1/2}_{\widetilde{\omega}}$ on $\h$ with 
 $\|V\|= \| S^{1/2}_{\widetilde{\omega}}\|$
 and
  $\|V^{-1}\|= \| S^{-1/2}_{\widetilde{\omega}}\|$.
 Then
 \begin{equation}\label{omeg1}
\widetilde{\omega_j} =
\sum_{i\in I} \< S_{\widetilde{f}}^{-1/2} \widetilde{f}_i, e_j\> V z_i, \ j\in I,
\end{equation}
 and furthermore,
 $\|V\|=
\sqrt{\|S_{\widetilde{\omega}}\|}\leq \sqrt{\|S_{\widetilde{f}}\|}$
and $\|V^{-1}\|=
 \sqrt{\|S_{\widetilde{\omega}}^{-1}\|}
 \leq\sqrt{\|S_{\widetilde{f}}^{-1}\|}$, implying that $\seqgrj[\widetilde{\omega}]$ is an R-dual of type III of $\{\widetilde{f_i}\}_{i\in I}$.

 Now assume in addition that $Q$ satisfies the property (\ref{eqstar}).
 By Theorem \ref{cordualopt}, $\otj$ has the same optimal bounds as $\fti$.
Then the optimal bounds of
$\seqgrj[\widetilde{\omega}]$ are the same as the optimal bounds of $\{\widetilde{f_i}\}_{i\in I}$, which by Proposition \ref{rdualopt} implies that the operator $V$ must satisfy a property analogue to  (\ref{eqstar}), precisely as stated in the proposition (note that the ranges of the analysis operators of  $\{\widetilde{f_i}\}_{i\in I}$ and $\fti$ are the same, and $S_{\widetilde{f}}=S^{-1}$).
\ep

\section{R-duals of Riesz bases}\label{rbcase}

In this section we consider the special case where $\fti$ is a Riesz basis for $\h.$  We will
provide extra insight into the various R-duals by
characterizing   the R-duals of type I, II, III, IV, as well as III having the essential property (\ref{eqstar}).

For the purpose of one of the characterizations,
   recall (see e.g. \cite{wigner}) that $G:\h\to\h$
 is called an {\it antiunitary transformation}\footnote{also called {\it antiunitary operator}  in the literature \cite{wigner}, but we do not use this term, because throughout the paper, the term  \lq\lq operator\rq\rq\,  is used with the meaning of a linear mapping}
 of the complex Hilbert space $\h$
  if it is a bijective 
  mapping of $\h$ onto $\h$ so that
 $\<Gx, Gy\>=\<y,x\>$ for all $x$ and $y$ in $\h$; such a transformation is as consequence antiliner, i.e.,  $G(ax+by)= \overline{a}Gx+ \overline{b}Gy $ for all $a,b\in\mc$ and all $x,y\in\h$.

\begin{prop} \label{pr1}
Let $\seqgri[f]$ be a Riesz basis for $\h$ with frame operator $S$ and optimal bounds $A,B$. Then the following statements hold.
\begin{itemize}

  \item[{\rm (i)}]
  The R-duals  of type I
of  $\seqgri[f]$ are precisely the Riesz bases $\seqgri[\omega]$   for which
there is an antiunitary transformation $G:\h\to\h$ so that
 $ \{S^{-1/2} G \omega_j\}_{j\in I}$ is an orthonormal basis for $\h;$
in particular, the optimal bounds of $\seqgri[\omega]$  are $A,B$.

\item[{\rm (ii)}]
 The R-duals  of type II
of  $\seqgri[f]$ are precisely the Riesz bases $\seqgri[\omega]$   for which
 $ \{S^{-1/2} \omega_j\}_{j\in I}$ is an orthonormal basis for $\h$ (implying that the optimal bounds of $\seqgri[\omega]$  are $A,B$).

  \item[{\rm (iii)}]
The R-duals  of type III of  $\seqgri[f]$ having the property (\ref{eqstar}) 
 are precisely the Riesz bases 
  for $\h$ with optimal bounds $A,B$.

 \item[{\rm (iv)}]
The R-duals  of type III
of  $\seqgri[f]$ are precisely the Riesz bases 
 which have $A, B$ as bounds.

  \item[{\rm (v)}]
The R-duals  of type IV of  $\seqgri[f]$
 are precisely the Riesz bases  for $\h$.

  \end{itemize}
\end{prop}

\bp (i) Let $\{z_i\}_{i\in I}$ be the orthonormal basis $\{S^{-1/2}f_i\}$.

First assume that $\seqgrj[\omega]$ is an R-dual of type I of $\seqgri[f]$ with respect to some orthonormal bases $\seqgri[e]$, $\seqgri[h]$. By Theorem \ref{207b}, $\seqgrj[\omega]$ is a Riesz sequence with optimal bounds $A, B$. 
 Consider the mapping $G$ determined by $Gh:=\sum_{i\in I}\< h_i, h\> z_i$, $h\in\h$. Then $G$ is an antiunitary transformation of $\h$ and  for every $j\in I$,
  \begin{eqnarray*}
  S^{-1/2}G \omega_j&=&
   S^{-1/2} \left(\sum_{i\in I} \<h_i,\omega_j\>z_i\right)
  =  S^{-1/2} \left(\sum_{i\in I} \<e_j, f_i\>z_i\right) \\ 
  &=&  S^{-1/2} \left(\sum_{i\in I} \<S^{1/2} e_j,  z_i\>z_i\right) 
  =e_j,
\end{eqnarray*}
 which leads to the desired conclusion.

Conversely, assume that $G:\h\to\h$ is an antiunitary transformation on $\h$ and $
 \{S^{-1/2} G\omega_j\}_{j\in I }
$
is an orthonormal basis of $\h$; denote this orthonormal basis by $\seqgrj[e]$. The sequence $\{G^{-1}z_i\}_{i\in I}$ is an orthonormal basis of $\h$ and the mapping $E$ given by $Eh:= \sum_{i\in I} \<h, z_i\> G^{-1}z_i$ is well defined from $\h$ into $\h$. Furthermore, observe that  $E$ is a unitary operator and $Ez_i=G^{-1}z_i$, $i\in I$. Define $h_i:= Ez_i$, $i\in I$. Then  $\{h_i\}_{i\in I}$ is an orthonormal basis of $\h$ and  for every $j\in I$,
\begin{eqnarray*}
\sum_{i\in I} \<f_i, e_j\>h_i
&=& \sum_{i\in I} \<f_i, S^{-1/2} G\omega_j\>E z_i
= \sum_{i\in I} \<S^{-1/2} f_i,  G\omega_j\>E z_i\\
&=& \sum_{i\in I} \<z_i, G\omega_j\> G^{-1}z_i
= G^{-1}(\sum_{i\in I} \<G\omega_j, z_i\> z_i)= \omega_j,
\end{eqnarray*}
which implies that $\seqgrj[\omega]$ is an R-dual of type I of $\seqgri[f]$.

(ii) First assume that $\seqgrj[\omega]$  is an R-dual of type II of $\seqgri[f]$. Then $\otj$ is a Riesz basis for $\h$ and by Lemma \ref{impcond}(iv),
 $ \{S^{-1/2} \omega_j\}_{j\in I}$ is an orthonormal basis for $\h$.
Now \cite[Lemma 1.2]{sc} implies that $\{\omega_j\}_{j\in I}$
has optimal bounds $A, B$.

For the converse, assume that $\seqgri[\omega]$ is a Riesz basis for $\h$ 
such that $ \{S^{-1/2} \omega_j\}_{j\in I}$ is an orthonormal basis for $\h$.
By Lemma \ref{impcond}(ii), $ \{S^{-1/2} \omega_j\}_{j\in I}$ is an R-dual of type I of $ \{S^{-1/2} f_i\}_{i\in I}$, i.e., there exist orthonormal bases $\seqgri[e]$ and $\seqgri[h]$ so that $S^{-1/2} \omega_j = \sumgrj \<S^{-1/2} f_j, e_i\>h_j$, $j\in I$, which implies that
$\omega_j = \sumgrj \<S^{-1/2} f_j, e_i\>S^{1/2}h_j$, $j\in I$.

(iii)  If $\seqgrj[\omega]$  is an R-dual of type III of $\seqgri[f]$ having
the property (\ref{eqstar}),
the conclusion follows from Theorem \ref{cordualopt}(iii).
Now assume that  $\seqgrj[\omega]$  is a Riesz basis for $\h$ having $A, B$ as optimal bounds. 
By \cite[Theorem 4.4(ii)]{sc} and the validity of (\ref{223f}), it follows that  $\seqgrj[\omega]$  is an R-dual of type III of $\seqgri[f]$ with respect to some triplet ($\seqgri[e], \seqgri[h], Q$). By Proposition \ref{rdualopt}, the property (\ref{eqstar}) must hold.

(iv) If $\seqgrj[\omega]$  is an R-dual of type III of $\seqgri[f]$, the conclusion follows from
\cite[Prop. 4.3]{sc}.
Now assume that  $\seqgrj[\omega]$  is a Riesz basis for $\h$ having $A, B$ as bounds. 
By \cite[Theorem 4.4(ii)]{sc} and the validity of (\ref{223f}), it follows that $\seqgrj[\omega]$  is an R-dual of type III of $\seqgri[f]$.

(v) If $\seqgrj[\omega]$  is an R-dual of type IV of $\seqgri[f]$, then $\seqgrj[\omega]$ is a Riesz basis for $\h.$
For the converse, assume that $\seqgri[\omega]$ is a Riesz basis for $\h$ and denote its frame operator by $S_\Omega$. 
By Lemma \ref{impcond}(ii), $ \{S_\Omega^{-1/2} \omega_j\}_{j\in I}$ is an R-dual of type I of $ \{S^{-1/2} f_i\}_{i\in I}$, i.e., there exist orthonormal bases $\seqgri[e]$ and $\seqgri[h]$ so that $S_\Omega^{-1/2} \omega_j = \sumgrj \<S^{-1/2} f_j, e_i\>h_j$, $j\in I$. Therefore,
$\omega_j = \sumgrj \< f_j, S^{-1/2}e_i\>S_\Omega^{1/2}h_j$, $j\in I$, which is an R-dual of type IV of $\fti$.
    \ep

 \vspace{.1in}
\noindent{\bf Acknowledgment:} Diana Stoeva  acknowledges support from the
 Austrian Science Fund (FWF) START-project FLAME ('Frames and Linear Operators for Acoustical Modeling and Parameter Estimation'; Y 551-N13).
She is grateful for the hospitality of the Technical University of Denmark, where much of the work on the paper was done.

\begin{tabbing}
text-text-text-text-text-text-text-text-text-text \= text \kill \\
Diana T. Stoeva  \> Ole Christensen \\
Acoustics Research Institute \> Technical University of Denmark\\
Wohllebengasse 12-14 \> DTU Compute\\
 Vienna A-1040, Austria\> 2800 Lyngby, Denmark  \\

Email: dstoeva@kfs.oeaw.ac.at \> Email: ochr@dtu.dk
\end{tabbing}

\end{document}